\documentclass[oneside,leqno,12pt]{article}
\usepackage{amssymb,amsmath,latexsym}

\def\gh{\mathfrak{h}}

\def\gn{\mathfrak{n}}
\def\go{\mathfrak{o}}
\def\gp{\mathfrak{p}}

\def\gs{\mathfrak{s}}

\def\gu{\mathfrak{u}}
\def\gv{\mathfrak{v}}
\def\gw{\mathfrak{w}}

\def\gz{\mathfrak{z}}
\def\gN{\mathfrak{N}}
\def\gZ{\mathfrak{Z}}

\def\Aut{{\rm Aut}}

\def\Ad{{\rm Ad}}
\def\ad{{\rm ad}}

\def\Ind{{\rm Ind\,}}

\def\Im{{\rm Im\,}}
\def\Re{{\rm Re\,}}

\def\Skew{{\rm Skew\,}}  % antisymmetric matrices
\def\C{\mathbb{C}}

\def\H{\mathbb{H}}

\def\O{\mathbb{O}}

\def\R{\mathbb{R}}

\def\cD{\mathcal{D}}

\def\cF{\mathcal{F}}

\def\cH{\mathcal{H}}

\def\cM{\mathcal{M}}

\def\cO{\mathcal{O}}
\def\cP{\mathcal{P}}

\def\cS{\mathcal{S}}

\newtheorem{theorem}[equation]{Theorem}

\newtheorem{corollary}[equation]{Corollary}

\begin{document}

\title{Classical Analysis and Nilpotent Lie Groups}

\author{Joseph A. Wolf} 

\date{}

\maketitle

Classical Fourier analysis has an exact counterpart in group theory
and in some areas of geometry.  Here I'll describe how this goes for
nilpotent Lie groups, for a class of Riemannian manifolds closely
related to a nilpotent Lie group structure.  There are also some 
infinite dimensional analogs but I won't go into that here.  
The analytic ideas are not so different from the classical
Fourier transform and Fourier inversion theories in one real variable.

Here I'll give a few brief indications of this beautiful topic.  References,
proofs and related topics for the finite dimensional theory can be found a 
recent AMS Monograph/Survey volume \cite{csbook}.
If you are interested in the infinite dimensional theory you may also wish 
to look at the article \cite{dl3} appear soon in Mathematische Annalen.

In Section \ref{sec1}
I'll recall a few basic facts on classical Fourier theory and note
the connection with the theory of locally compact abelian groups
and their unitary representations.
In Section \ref{sec2} we look at the first noncommutative locally
compact groups, the Heisenberg groups $H_n$.  We describe their unitary
representations, Fourier transform  theory and Fourier inversion
formula.  

The coadjoint orbit picture is the best way to understand
representations of nilpotent Lie groups.  It is guided by the
example of the Heisenberg group.  We indicate that theory in
Section \ref{sec3}.
Then in Section \ref{sec4} we come to a class of connected,
simply connected, nilpotent Lie groups with many of the good analytic
properties of vector groups and Heisenberg groups.  Those are the
simply connected, nilpotent Lie groups with square integrable
representations.

In Section \ref{sec5} we push some of the analysis to a class of
homogeneous spaces where the techniques and results are analogous to
those of locally compact abelian groups.  Those are the commutative
spaces $G/K$, i.e. the Gelfand pairs $(G,K)$.  We have already had a
glimpse of this in Section \ref{sec2} for the semidirect products
$H_n \rtimes K$ and the riemannian homogeneous spaces $(H_n \rtimes K)/K$
where $K$ is a compact group of automorphisms of $H_n$.  In Section
\ref{sec6}  we look
more generally at Fourier transform  theory and the Fourier inversion
formula for commutative nilmanifolds $(N\rtimes K)/K$ where
$N$ is a simply connected, nilpotent Lie groups with square integrable
representations and $K$ is a compact group of automorphisms of $N$.

As indicated earlier, Fourier analysis for
nilpotent Lie groups $N$ and commutative nilmanifolds $(N\rtimes K)/K$, 
where $N$ has square integrable representations, has recently been
extended to some classes direct limit groups and spaces.  

\section{Classical Fourier Series}\label{sec1}

Let's recall the Fourier series development for a function $f$ of one 
variable that is periodic of period $2\pi$. One views $f$ 
as a function on the circle  $S = \{e^{ix}\}$.
The circle $S$ is a multiplicative group and we expand $f$
in terms of the unitary characters
$$
\chi_n: S \to S \text{ by } \chi_n(x) = e^{i nx}\,\, ,
\text{ continuous group homomorphism.}
$$ 
Then the Fourier inversion formula is
$$
f(x) = \sum_{n = -\infty}^\infty \widehat{f}(n) \chi_n
$$
where the Fourier transform
$$
\widehat{f}(n) = \tfrac{1}{2\pi} \int_0^{2\pi} f(x)e^{-i nx} dx
= \langle f, \chi_n \rangle_{L^2(S)}.
$$
The point is that $f$ is a linear combination of the $\chi_n$ with
coefficients given by the Fourier transform $\widehat{f}$.  This uses
the topological group structure and the rotation--invariant measure
on $S$.
\medskip

One has a similar situation when the compact group $S$ is replaced by a
finite dimensional real vector space $V$.  Let $V^*$ denote its linear
dual space.  If $f \in L^1(V)\cap L^2(V)$ the Fourier inversion formula is
$$
f(x) =  \left (\tfrac{1}{2\pi}\right )^{m/2}\int_{V^*} 
\widehat{f}(\xi) e^{ix\cdot \xi} d\xi
$$
where the Fourier transform is
$$
\widehat{f}(\xi) = \left (\tfrac{1}{2\pi}\right )^{m/2}
\int_V f(x)e^{-i x\cdot \xi} dx = \langle f, \chi_\xi \rangle_{L^2(V)}.
$$
Again, $f$ is a linear combination (this time it is a continuous linear
combination) of the unitary characters $\chi_\xi(x) = e^{ x\cdot \xi}$
on $V$, and the coefficients of the linear combination are given by the
Fourier transform $\widehat{f}$.
\medskip

It is the same story for locally compact abelian groups $G$. 
The unitary characters form a group
$$
\widehat{G} = \{\chi: G \to S \text{ continuous homomorphisms}\}
$$
with composition $(\chi_1\chi_2)(x) = \chi_1(x)\chi_2(x)$.  It is
locally compact with the weak topology for the evaluation maps 
$ev_x: \chi \mapsto \chi(x)$.
If $f \in L^1(G)\cap L^2(G)$ the Fourier inversion formula is
$$
f(x) = \int_{\widehat{G}} \widehat{f}(\chi) \chi(x) d\chi
$$
where the Fourier transform
$$
\widehat{f}(\chi) = \int_G f(x)\overline{\chi(x)}dx
 = \langle f, \chi\rangle_{L^2(G)}.
$$
As in the Euclidean cases, Fourier inversion expresses the function $f$ as a
(possibly continuous) linear combination of unitary characters on
$G$, where the coefficients of the linear combination are given by the
Fourier transform.
\medskip

In this context, $f \mapsto \widehat{f}$ preserves $L^2$ norm and extends 
by continuity to an isometry of $L^2(G)$ onto $L^2(\widehat{G})$.  In effect
this expresses $L^2(G)$ as a (possibly continuous) sum of $G$--modules,
$$
L^2(G) = \int_{\widehat{G}} \C_\chi d\chi \text{ where } \C_\chi 
\text{ is spanned by } \chi.
$$
In this direct integral decomposition $d\chi$ could be replaced by any
equivalent measure, so that decomposition is
not as precise as the Fourier inversion formula.

\section{The Heisenberg group}\label{sec2}

Next, we see what happens when we weaken the commutativity condition.
The first case of that is the case of the Heisenberg group.  
There the Fourier transform and Fourier inversion are in some sense the
same as in the classical case of a vector group, except that some of the
integration occurs in the character formula and the rest in integration
over the unitary dual.

The Heisenberg group of real dimension $2m+1$ is
$$
H_m = \Im\C + \C^m \text{ with group law }
(z,w)(z',w') = (z+z'+\Im \langle w,w'\rangle, w+w')
$$
where $\Im$ denotes imaginary component (as opposed to the coefficient
of $\sqrt{-1}$), $z,z' \in Z := \Im\C$ and $w,w' \in W := \C^m)$. 
Its Lie algebra, the Heisenberg algebra, 
is 
$$
\gh_m = \gz + \gw = \Im\C + \C^m \text{ with }
[z+w,z'+w'] = (z+z'+2\,\Im\langle w,w'\rangle) + (w+w').
$$
Here $Z = \exp(\gz)$ is both the center and the derived group of
$H_m$, and its complement $W = \exp(\gw) \cong \R^{2m}$.
\medskip

Unitary characters have to annihilate the derived group of $H_m$, in
other words factor through $H_m/Z$, so the only functions that can
be expanded in unitary characters are the functions that are lifted
from $H_m/Z$.  Thus we have to consider something more general.
\medskip

The space $\widehat{H_m}$ of (equivalence classes of) irreducible unitary
representations of $H_m$ breaks into two parts, one consisting of the
$1$--dimensional representations and the other of the infinite dimensional
representation.  This goes as follows.
\begin{itemize}
\item One-dimensional representations.  These are the ones that annihilate the
center $Z$, and are given by the unitary characters $\chi_\xi$,
$\xi \in W^*$, on $W \cong \R^{2m}$.
\item Infinite dimensional representations.  These are the 
$\pi_\zeta = \Ind_N^{H_m}(\chi_\zeta)$ where
$$
N = \Im\C + \R^m \subset H_m \text{ and } \zeta \in \gz^*\setminus \{0\}\, .
$$
\end{itemize}
Recall the definition of the induced representation 
$\pi_\zeta = \Ind_N^{H_m}(\chi_\zeta)$.
Here $\chi_\zeta$ extends from $Z$ to $N$ by 
$\chi_\zeta(z,w) = \chi_\zeta(z)$.  Thus we have a unitary line bundle
over $H_m/N$ associated to the principal $N$--bundle
$H_m \to H_m/N$ by the action $w: t \mapsto \chi_\zeta(w)t$ of
$N$ on $\C$.  Now $\pi_\zeta$ is the natural action of $H_m$  on the 
space of $L^2$ sections of that line bundle.

The classical
``Uniqueness of the Heisenberg commutation relations'' says that $\zeta$
determines the class $[\pi_\zeta] \in \widehat{H_m}$.
And restriction to $Z$ shows that $[\pi_\zeta] = [\pi_{\zeta'}]$
just when $\zeta = \zeta'$.

Using the fact that $\zeta$ determines $[\pi_\zeta]$, one realizes 
$[\pi_\zeta]$ by an action of $H_m$ on the Hilbert space $\cH_m$ of Hermite
polynomials on $\C^m$.  The maximal compact subgroup of $\Aut(H_m)$ is the
unitary group $U(m)$.  Its action is
$$
g : (z,w) \mapsto (z,g(w)).
$$
Result: $\pi_\zeta$ extends to an irreducible unitary representation
$\widetilde{\pi_\zeta}$ of the semidirect product $H_m\rtimes U(m)$ on $\cH_m$.
So if $K$ is any closed subgroup of $U(m)$ then
$\widetilde{\pi_\zeta}|_{H_m\rtimes K}$ is an irreducible unitary
representation of $H_m\rtimes K$ on $\cH_m$.  The Mackey Little Group
theory says that $\widehat{H_m\rtimes K} = \{[\widetilde{\pi_\zeta} \otimes
\kappa] \mid [\kappa] \in \widehat{K} \text{ and } 
\zeta \in \gz^*\setminus \{0\}\}$.
 
\section{Representations and Coadjoint Orbits}\label{sec3}

Kirillov theory for connected simply connected nilpotent Lie groups $N$
realizes their unitary representations in terms of the the coadjoint
representation of $N$, that is, the representation $\ad^*$ of $N$
on the linear dual space $\gn^*$ of its Lie algebra $\gn$.

On the group level the coadjoint representation is given by $(\Ad^*(n)f)(\xi) 
= f(\Ad(n)^{-1}\xi)$.  Write $\cO_f$ for the (coadjoint) orbit $\Ad^*(N)f$
of the linear functional $f$.  Consider the antisymmetric bilinear form
$b_f$ on $\gn^*$ given by $b_f(\xi,\eta) = f([\xi,\eta])$. The kernel 
of $b_f$ is the Lie algebra of the isotropy subgroup of $N$ at $f$.  Thus
$b_f$ defines an $\Ad^*(N)$--invariant symplectic form $\omega_f$ on the
coadjoint orbit $\cO_f$.  The symplectic homogeneous space $(\cO_f,\omega_f)$ 
leads to a 
unitary representation class $[\pi_f] \in \widehat{N}$, as follows.

Let $N_f$ denote the $\Ad^*(N)$--stabilizer of $f$.  Its Lie algebra $\gn_f$
is the annihilator of $f$, in other words 
$\gn_f = \{\nu \in \gn \mid f(\nu,\gn) = 0\}$.  A (real) {\bf polarization}
for $f$ is a subalgebra $\gp \subset \gn$ that contains $\gn_f$, has
dimension given by $\dim(\gp/\gn_f) = \tfrac{1}{2}\dim(\gn/\gn_f)$, and
satisfies $f([\gp,\gp]) = 0$.  Under the differential $\gn \to T_f(\cO_f)$
of $N \to \cO_f$, real polarizations for $f$ are in
one to one correspondence with $N$--invariant integrable Lagrangian
distributions on $(\cO_f,\omega_f)$.

Fix a real polarization $\gp$ for $f$ and let $P = \exp(\gp)$.  It is the
analytic subgroup of $N$ for $\gp$ and it is a closed, connected, simply
connected subgroup of $N$.  In particular $e^{if}: P \to \C$ is a well
defined unitary character.  That defines a unitary representation
$$
\pi_f = \pi_{f,\gp} = \Ind_P^N(e^{if})
$$
of $N$.  The basic facts are given by
\begin{theorem}\label{kir-thm} 
Let $N$ be a connected simply connected nilpotent Lie group and $f \in \gn^*$.

{\rm 1.} There exist real polarizations $\gp$ for $f$.

{\rm 2.} If $\gp$ is a real polarization for $f$ then the unitary 
representation $\pi_{f,\gp}$ is irreducible.

{\rm 3.} If $\gp$ and $\gp'$ are real polarizations for $f$ then the unitary
representations $\pi_{f,\gp}$ and $\pi_{f,\gp'}$ are equivalent, so the
class $[\pi_f] \in \widehat{N}$ is well defined.

{\rm 4.} If $[\pi] \in \widehat{N}$ then there exists $h \in \gn^*$ such
that $[\pi] = [\pi_h]$.
\end{theorem}
In other words, $f \mapsto \pi_{f,\gp}$ induces a one to one map of
$\gn^*/\Ad^*(N)$ onto $\widehat{N}$.

To see just how this works, consider the case where $N$ is the Heisenberg
group $H_m$, and let $f \in \gh_m^*$.  Here the center
$\gz = \Im \C$ and its complement $\gv = \C^n$.  Decompose 
$\gv = \gu + \gw$ where $\gu = \R^n$ and $\gw = i\R^n$.  Note 
$\Im \langle \gu,\gu \rangle = 0 = \Im \langle \gw,\gw \rangle$.
If $f(\gz) = 0$ we have the real polarization 
$\gp = \gh_m$.  If $f(\gz) \ne 0$ we have the real polarization $\gp = 
\gz + \gu$.  That demonstrates Theorem \ref{kir-thm}(1).  

If $f(\gz) = 0$ then $\pi_{f,\gp}$ is a unitary character on $H_n$, 
automatically irreducible.  Now suppose $f(\gz) \ne 0$ and $\gp =
\gz + \gu$.  Then $\pi_{f,\gp}$ is a representation of $H_n$ on
$L^2(G/P) = L^2(W)$ where $W = \exp(\gw) \cong \R^n$.  Then $\pi_{f,\gp}(H_n)$
acts by all translations on $W$ and by scaling that distinguishes the
integrands of $L^2(W) = \int_{\gw^*}\, e^{\langle \xi, \cdot \rangle}\C\, d\xi$,
so it is irreducible.  That demonstrates Theorem \ref{kir-thm}(2).

If $f(\gz) = 0$ then $\gh_n$ is the only real polarization for $f$.  Now
suppose $f(\gz) \ne 0$ and consider the case where $\gp = \gz + \gu$ and
$\gp' = \gz + \gw$.  Then $\omega = \Im h(\cdot,\cdot)$ pairs $\gu$ with
$\gw$, and the Fourier transform $\cF: L^2(U) \cong L^2(W)$ intertwines
$\pi_{f,\gp}$ with $\pi_{f',\gp}$.  More generally, if $\gp$ and $\gp'$
are any two real polarizations for $f$, then we write $\gp = (\gp\cap\gp') + 
\gu'$, $\gp' = (\gp\cap\gp') + \gw'$ and $\gp\cap\gp' = \gz + \gv'$.  That 
done, $\omega$ pairs $\gu'$ with $\gw'$, and the corresponding Fourier 
transform $\cF: L^2(U') \cong L^2(W')$ combines with the identity
transformation of $L^2(V')$ to give a map
$L^2(V')\widehat{\otimes} L^2(U') \cong L^2(V')\widehat{\otimes} L^2(W')$
that intertwines $\pi_{f,\gp}$ with $\pi_{f',\gp}$.  This demonstrates
Theorem \ref{kir-thm}(3).

Now Theorem \ref{kir-thm}(4) follows from the considerations we
outlined in Section \ref{sec2}.
In the terminology there,
the infinite dimensional irreducible unitary representation $\pi_\zeta$ 
of $H_n$ is equivalent to $\pi_f$ whenever $f \in \gn^*$ such that
$f(z) = \langle\zeta, x\rangle$ for every $z\in \gz$.  In particular, if
$f(\gz) \ne 0$ where $\gz$ is the center of the Heisenberg algebra, then
the coadjoint orbit $\cO_f = f + \gz^\perp$, where
$\gz^\perp := \{h \in \gn^* \mid h(\gz) = 0\}$.
Of course one can also verify this by direct computation.

\section{Square Integrable Representations}\label{sec4}

In this section $N$ is a connected simply connected nilpotent Lie group
and $Z$ is its center.  
If $\zeta \in \widehat{Z}$ we denote
$\label{rel-dual}
\widehat{N}_\zeta = \{[\pi] \in \widehat{N} \mid \pi|_Z
	\text{ is a multiple of } \zeta$.
The corresponding $L^2$ space is
\begin{equation}\label{rel-l2}
% \begin{aligned}
L^2(N/Z:\zeta) := 
\left \{f: N \to \C \text{ measurable } \left |\, 
	\begin{matrix} f(nz) = \zeta(z)^{-1}f(n) \text{ and }\\
	\phantom{X} \int_{N/Z} |f(n)|^2 d\mu_{_{N/Z}}(nZ) < \infty
	\end{matrix} \right . \right \}.
% \end{aligned}
\end{equation}
The inner product 
$\langle f, h \rangle_\zeta = \int_{N/Z} f(n)\overline{h(n)} d\mu_{_{N/Z}}(nZ)$
is well defined on the relative $L^2$ space $L^2(N/Z:\zeta)$.
Each $\widehat{N}_\zeta$ is a measurable subset of $\widehat{N}$, and 
$\widehat{N} = \bigcup_{\zeta \in \widehat{Z}}\, \widehat{N}_\zeta$.
Here $L^2(N) = \int_{\widehat{Z}} L^2(N/Z:\zeta)$.  This decomposes
the left regular representation of $N$ as 
$$
\ell = \Ind_{\{1\}}^N(1) = \Ind_Z^N \Ind_{\{1\}}^Z(1) 
= \Ind_Z^N \int_{\widehat{Z}}\,\zeta \, d\zeta
= \int_{\widehat{Z}}\Ind_Z^N\zeta \, d\zeta
= \int_{\widehat{Z}}\, \ell_\zeta \, d\zeta
$$ 
where $\ell_\zeta = \Ind_Z^N\zeta$ is the
left regular representation of $N$ on $L^2(N/Z:\zeta)$.  The corresponding
expansion for functions, 
$$
f(n) = \int_{\widehat{Z}} f_\zeta(n) d\zeta 
\text{ where } f_\zeta(n) = \int_Z f(nz)\zeta(z) d\mu_{_Z}(z),
$$ 
is just Fourier inversion on the commutative locally compact group $Z$.

Now we describe some results of Moore and myself on square integrable
representations in this context.  The first observation is

\begin{theorem}\label{mw-A} 
Let $N$ be a connected simply connected nilpotent Lie group and
$\zeta \in \widehat{Z}$.  If $[\pi] \in \widehat{N}_\zeta$ then 
the following conditions are equivalent.

{\bf 1.}  There exist nonzero $u,v \in H_\pi$ such that $|f_{u,v}|
\in L^2(N/Z)$, i.e., $f_{u,v} \in L^2(N/Z:\zeta)$.

{\bf 2.}  The coefficient $|f_{u,v}| \in L^2(N/Z)$, equivalently
$f_{u,v} \in L^2(N/Z:\zeta)$, for all $u,v \in H_\pi$.

{\bf 3.}  $[\pi]$ is a discrete summand of $\ell_\zeta$.
\end{theorem}

A representation class $[\pi] \in \widehat{N}$ is
$\text{\bf {\sl L}}^2$ or {\bf square integrable} or {\bf relative
discrete series} if its coefficients $f_{u,v}(n) = f_{\pi:u,v}(n) 
:= \langle u, \pi(n)v\rangle$ satisfy $|f_{u,v}| \in L^2(N/Z)$,
in other words if its coefficients are square integrable modulo $Z$.  
Theorem \ref{mw-A} says that it is sufficient to check this for just
one nonzero coefficient, and Theorem \ref{mw-A}(3) justifies the term
``relative discrete series''.  

We say that $N$ {\bf has square integrable representations} if at least one
class $[\pi] \in \widehat{N}$ is square integrable.  
These representations satisfy an analog of the Schur
orthogonality relations:

\begin{theorem}\label{mw-B} 
Let $N$ be a connected simply connected nilpotent Lie group. If
$\zeta \in \widehat{Z}$ and $[\pi] \in \widehat{N}_\zeta$ is
square integrable then there is a number $\deg(\pi) > 0$ 
such that the coefficients of $\pi$ satisfy
\begin{equation}
\int_{N/Z} f_{u_1,v_1}(n)\overline{f_{u_2,v_2}(n)} d\mu_{_{N/Z}}(nZ) =
	\tfrac{1}{\deg(\pi)}\langle u_1,u_2\rangle \,
	\overline{\langle v_1,v_2\rangle}
\end{equation}
for all $u_i, v_i \in H_\pi$.  If $[\pi_1], [\pi_2] \in \widehat{N}_\zeta$
are inequivalent square integrable representations then their coefficients
are orthogonal in $L^2(N/Z:\zeta)$,
\begin{equation}
\int_{N/Z} f_{\pi_1: u_1,v_1}(n)\overline{f_{\pi_2: u_2,v_2}(n)} 
	d\mu_{_{N/Z}}(nZ) = 0,
\end{equation}
for all $u_1, v_1 \in H_{\pi_1}$ and $u_2, v_2 \in H_{\pi_2}$.
\end{theorem}

The number $\deg(\pi)$ is the {\bf formal degree} of $[\pi]$.  It plays
the same role in Theorem \ref{mw-B} as that played by the degree in the
Schur orthogonality relations for compact groups .  In general,
$\deg(\pi)$ depends on normalization of Haar measure: a rescaling of
Haar measure $\mu_{_{N/Z}}$ of $N/Z$ to $c\mu_{_{N/Z}}$ rescales formal
degrees $\deg(\pi)$ to $\tfrac{1}{c}\deg(\pi)$.  We don't see this for
compact groups because there we always scale Haar measure to total mass $1$.

Theorems \ref{mw-A} and \ref{mw-B} only require that $N$ be a locally compact 
group of Type I and that $Z$ be a closed subgroup of the center of $N$.  They 
can be understood as special cases of Hilbert algebra theory.
Here we related them to the Kirillov theory.  

Given $f \in \gn^*$ we have the bilinear form $b_f(x,y) = f([x,y])$, 
the coadjoint orbit
$\cO_f = \Ad^*(N)f$, the associated representation $[\pi_f]$, and the
character $\zeta \in \widehat{Z}$ such that $[\pi_f] \in \widehat{N}_\zeta$.
Note that $f|_{\gz}$ determines the affine subspace $f + \gz^\perp$ in $\gn^*$.

\begin{theorem}\label{mw-1} 
Let $N$ be a connected simply connected nilpotent Lie group and $f \in \gn^*$.
Then the following conditions are equivalent.

{\bf 1.} $[\pi_f]$ is square integrable.

{\bf 2.} The left regular representation $\ell_\zeta$ of $N$ on 
$L^2(N/Z:\zeta)$ is primary.

{\bf 3.} $\cO_f = f + \gz^\perp$, determined by the restriction $f|_\gz$.

{\bf 4.} $b_f$ is nondegenerate on $\gn/\gz$.
\end{theorem}

Recall the notion of the Pfaffian $\text{\rm Pf}(\omega)$ 
of an antisymmetric bilinear form
$\omega$ on a finite dimensional real vector space $V$ relative to a
volume form $\nu$ on $V$. If $\dim V$ is odd then by definition 
$\text{\rm Pf}(\omega) = 0$.  If $\dim V = 2m$ even, then $\omega^m$ is a
multiple of $\nu$, and by definition that multiple of $\text{\rm Pf}(\omega)$;
in other words $\omega^m = \text{\rm Pf}(\omega)\nu$.  The Pfaffian is the
square root of the determinant on antisymmetric bilinear forms.

Fix a volume element $\nu$ on $\gv := \gn/\gz$.
If $f \in \gn^*$ we view $\omega_f(x,y) = f([x,y])$ as an antisymmetric
bilinear form on $\gv$.  Define $P(f) := \text{\rm Pf}(\omega_f)$.
Then $P$ is a homogeneous polynomial function on $\gn^*$, and $P(f)$ 
depends only on $f|_\gz$a  So there is a homogeneous polynomial
function (which we also denote $P$) 
on $\gz^*$ such that $P(f) = P(f|_\gz)$.

In the case of the Heisenberg group $H_m$, $P$ is the homogeneous polynomial 
$P(\zeta) = \zeta(z_0)^m$ of degree $m$ on $\gz^*$.  Here the choice of
nonzero $z_0 \in \gz$ is a normalization, in effect a choice of unit vector.
As described in Theorem \ref{mw-4} below, this also gives the formal degree
of $[\pi_f]$ where $\zeta = f|_\gz$.  Further, as described in
Theorem \ref{mw-6}, it gives the Plancherel measure on $\widehat{H_m}$.

In view of Theorem \ref{mw-1} we now have

\begin{theorem}\label{mw-2} 
The representation $\pi_f$ is square integrable if and only if 
the Pfaffian polynomial $P(f|_\gz) \ne 0$.  In particular 
$\phi: f|_\gz \mapsto [\pi_f]$
defines a bijection from $\{\lambda \in \gz^* \mid P(\lambda) \ne 0\}$
onto $\{[\pi] \in \widehat{N} \mid [\pi] \text{ is square integrable}\}.
$
\end{theorem}

One can view the polynomial $P$ as an element of the symmetric
algebra $S(\gz)$, and since $\gz$ is commutative that symmetric algebra
is the same as the universal enveloping algebra $\gZ$.  From this one
can prove

\begin{corollary}\label{mw-3} 
The group $N$ has square integrable representations if and only if
the inclusion $\gZ \hookrightarrow \gN$ of universal enveloping 
algebras, induced by $\gz \hookrightarrow \gn$,
maps $\gZ$ onto the center of $\gN$.
\end{corollary}

Both formal degree and the polynomial $P$ are scaled by $1/c$ when Haar
measure on $N/Z$ is scaled by $c$, so the following is independent of
normalization of Haar measure on $N/Z$.

\begin{theorem}\label{mw-4} 
The formal degree of a square integrable representation 
$[\pi_f] = \phi(f|_\gz)$ is given by $\deg(\pi_f) = |P(f|_\gz)|$.
\end{theorem}

As in the semisimple case, the {\bf infinitesimal character} of a
representation class $[\pi] \in \widehat{N}$ is the associative algebra
homomorphism $\xi_\pi : Cent(\gN) \to \C$ from the center of the
enveloping algebra, such that $d\pi(\zeta)$ is scalar multiplication
by $\chi_\pi(\zeta)$.  (Initially this holds only on $C^\infty$ vectors,
but they are dense in $H_\pi$, so by continuity it holds on all vectors.)
If $\zeta \in \gz$ then $\chi_{\pi_f}(\zeta) = if(\zeta)$.  Now, from
Theorem \ref{mw-4},

\begin{corollary}\label{mw-5} 
If $[\pi] \in \widehat{N}$ then the formal degree $\deg(\pi) = |\chi_\pi(P)|$
where we understand the formal degree of a non square integrable representation
to be zero.
\end{corollary}

For the Plancherel formula and Fourier inversion we must normalize Haar 
measures.  Choose Haar measures $\mu_{_Z}$ and $\mu_{_{N/Z}}$; they define 
a Haar measure $\mu_{_N}$ by $d\mu_{_N} = d\mu_{_{N/Z}}\, d\mu_{_Z}$, i.e.
$$
\int_N f(n) d\mu_{_N}(n) = 
\int_{N/Z} \left ( \int_Z f(nz) d\mu_{_Z}(z)\right ) d\mu_{_{N/Z}}(nZ).
$$
Now we have Lebesgue measures $\nu_{_Z}$, $\nu_{_{N/Z}}$ and $\nu_{_N}$
on $\gz$, $\gn/\gz$ and $\gn$ specified by the condition that the 
exponential map have Jacobian $1$ at $0$, and they satisfy
$d\nu_{_N} = d\nu_{_{N/Z}}\, d\nu_{_Z}$.  Normalize Lebesgue measures on
the dual spaces by the condition that Fourier transform is an isometry;
that gives Lebesgue measures $\nu^*_{_Z}$, $\nu^*_{_{N/Z}}$ and $\nu^*_{_N}$
such that $d\nu^*_{_N} = d\nu^*_{_{N/Z}}\, d\nu^*_{_Z}$.

\begin{theorem}\label{mw-6} 
Let $N$ have square integrable representations.  Let
$c = m!2^m$ where $2m$ is the maximum dimension of the
$\Ad^*(N)$--orbits in $\gn^*$.  Then Plancherel measure for $N$ is
concentrated on the square integrable classes
and its image under the map 
$$
\phi^{-1}: 
\{[\pi] \in \widehat{N} \mid [\pi] \text{ is square integrable}\}
\to \{\lambda \in \gz^* \mid P(\lambda) \ne 0\}
$$
of {\rm Theorem \ref{mw-2}} is $c|P(x)|d\nu^*_{_Z}(x)$.
\end{theorem}

\section{Commutative Spaces -- Generalities}\label{sec5}

We just described the Fourier transform and Fourier inversion formulae for
$H_m$ --- and a somewhat larger class of connected simply connected nilpotent 
Lie groups.  Now we edge toward a more geometric setting, that of commutative
spaces, which is a common
generalization of Riemannian symmetric spaces, locally compact abelian groups
and homogeneous graphs.  It is interesting and precise for the cases that
involve connected simply connected nilpotent Lie groups with square 
integrable representations.

A {\it commutative space} $G/K$, equivalently a {\it Gelfand pair} $(G,K)$,
consists of a locally compact group $G$ and a compact subgroup $K$ such that
the convolution algebra $L^1(K\backslash G/K)$ is commutative.  There are
several other formulations.  Specifically, the following are equivalent.

1. $(G,K)$ is a Gelfand pair, i.e. $L^1(K\backslash G/K)$ is
commutative under convolution.

2. If $g,g' \in G$ then $\mu_{_{KgK}} * \mu_{_{Kg'K}} =
 \mu_{_{Kg'K}} * \mu_{_{KgK}}$
{\rm (convolution of measures on $K\backslash G/K$).}

3. $C_c(K\backslash G/K)$ is commutative under convolution.

4. The measure algebra $\cM(K\backslash G/K)$ is commutative.

5. The (left regular) representation of $G$ on $L^2(G/K)$ is multiplicity free.
\medskip

\noindent
If $G$ is a connected Lie group one can add

6. The algebra $\cD(G,K)$ of $G$--invariant
differential operators on $G/K$ is commutative.  
\medskip

Commutative spaces $G/K$ are important for a number of reasons.  First,
they are manageable because their basic harmonic analysis is very similar 
to that of locally compact abelian groups.  We will describe that in a
moment.  Second, in the Lie group cases, most of the $G/K$ carry invariant 
weakly symmetric Riemannian metrics, which have properties very similar to 
those of Riemannian symmetric spaces.  Third, the invariant differential
operators and corresponding spherical functions play a definite role
in special function theory.  And fourth, in the nilpotent Lie groups setting,
there is some interesting interplay between geometry and hypoellipticity.

We only consider the basic harmonic analysis, here in general and in
Section \ref{sec6} for the case of commutative nilmanifolds.

Analysis on locally compact abelian groups is based on decomposition of
functions in terms of unitary characters.  In the classical euclidean
case these are just the complex exponentials $\chi_\xi : V \to \C$,
$\xi \in V^*$, given by $\chi_\xi = e^{ix\cdot\xi}$. For a commutative space 
$G/K$ the appropriate replacements are the
positive definite spherical functions, defined as follows.  

A continuous 
$K$--bi--invariant function $\varphi: G \to \C$ is $K$--{\it spherical}
if $\varphi(1) = 1$ and $f \mapsto (f*\varphi)(1)$ is a homomorphism 
$C_c(K\backslash G/K) \to \C$.  Equivalent: $\varphi$ is not identically zero,
 and if $g_1, g_2 \in G$ then 
$$
\varphi(g_1)\varphi(g_2) = \int_K \varphi(g_1kg_2)dk.
$$
A function
$\varphi: G \to \C$ is {\it positive definite} if 
$\sum \varphi(g_j^{-1}g_i)\overline{c_j}c_i \geqq 0$ whenever 
$\{c_1, \dots , c_n\} \subset \C$ and $\{g_1, \dots , g_n\} \subset G$.

Denote $\cP = \cP(G,K)$: positive definite $K$--spherical functions on $G$.
There is a one--one correspondence $\varphi \leftrightarrow \pi_\varphi$
between $\cP$ and the irreducible unitary representations $\pi$ of $G$ that 
have a $K$--fixed unit vector $v$.  It is given by
$\varphi(g) = \langle v, \pi(g)v\rangle_{\cH_\pi}$.  We have the {\it spherical
transform} 
$$
\cS:f \mapsto \widehat{f} \text{ from } L^1(K\backslash G/K)
\text{ to functions on } \cP
$$
defined by
$$
\cS(f)(\varphi) = \widehat{f}(\varphi) = (f*\varphi)(1) = 
\int_G f(g)\varphi(g^{-1})dg.
$$
The corresponding spherical inversion formula is
$$
f(g) = \int_\cP \widehat{f}(\varphi) \varphi(g) d\mu(\varphi).
$$
Here $\cP$ has natural structure of locally compact space and $\mu$ is
called Plancherel measure. The spherical transform
$$
\cS: L^1(K\backslash G/K)\cap L^2(K\backslash G/K) 
\to L^2(\cP,\mu)
$$
preserves $L^2$ norm and extends by continuity to an isometry 
$$
\cS: L^2(K\backslash G/K)\cong L^2(\cP,\mu).
$$
Note that $\cS$ can only be given by its defining integral expression
when that integral converges.  This is why it has to be extended by
$L^2$ continuity.  Of course this problem is already present with the classical
Fourier transform on $\R$.

The Plancherel Formula $\cS: L^2(K\backslash G/K) \cong L^2(\cP,\mu)$ gives 
a continuous direct sum (direct integral) decomposition
$$
L^2(K\backslash G/K) \cong \int_\cP \C\varphi\,\, d\mu(\varphi).
$$
This extends to a continuous direct sum decomposition
$$
L^2(G/K) \cong \int_\cP \cH_{\pi_\varphi}\,\, d\mu(\varphi).
$$
Of course all this depends on knowledge of the Plancherel measure $\mu$.

\section{Commutative Nilmanifolds}\label{sec6}

Theorem of Carcano (special case): Let $K \subset U(m)$ acting on $\C^m$, 
where $\gh_m = \Im\C + \C^m$ with center $\Im\C$.
Then $(H_m\rtimes K)/K$ is commutative if and only if
the representation of $K_{_\C}$, on the ring of all polynomials on $\C^m$, is
multiplicity free.
\medskip

Ka\v c classified the connected $K_m$ that are irreducible on $\C^m$:
\medskip

{\footnotesize
\centerline{
\begin{tabular}{|c|l|l|}\hline
 & Group $K_m$ & Acting on  \\ \hline
1 & $SU(m)$ & $\C^m$,  $m \geqq 2$ \\ \hline
2 & $U(m)$ & $\C^m$, $m \geqq 1$ \\ \hline
3 & $Sp(n)$ & $\C^m = \C^{2n}$ \\ \hline
4 & $U(1) \times Sp(n)$ & $\C^m = \C^{2n}$ \\ \hline
5a & $U(1) \times SO(2n)$ & $\C^m = \C^{2n}$ \\ \hline
5b & $U(1) \times SO(2n+1)$ & $\C^m = \C^{2n+1}$ \\ \hline
6 & $U(n), n \geqq 2$ & $\C^m = S^2(\C^n)$  \\ \hline
7 & $SU(n), n$ odd & $\C^m = \Lambda^2(\C^n)$  \\ \hline
\end{tabular}
\begin{tabular}{|c|l|l|}\hline
 & Group $K_m$ & Acting on  \\ \hline
8 & $U(n)$ & $\C^m = \Lambda^2(\C^n)$ \\ \hline
9 & $SU(\ell) \times SU(n)$ &
        $\C^m = \C^\ell \otimes \C^n$, $\ell \ne n$ \\ \hline
10 & $S(U(\ell) \times U(n))$ &
        $\C^m = \C^\ell \otimes \C^n$  \\ \hline
11 & $U(2) \times Sp(n)$ &
        $\C^m = \C^2 \otimes \C^{2n}$ \\ \hline
12 & $SU(3) \times Sp(n)$ &
        $\C^m = \C^3 \otimes \C^{2n}$ \\ \hline
13 & $U(3) \times Sp(n)$ & $\C^m = \C^3 \otimes \C^{2n}$ \\ \hline
14 & $SU(n) \times Sp(4)$ & $\C^m = \C^n \otimes \C^8$ \\ \hline
15 & $U(n) \times Sp(4)$ & $\C^m = \C^n \otimes \C^8$ \\ \hline
\end{tabular}
}
}

A {\it commutative nilmanifold} is a commutative space $G/K$ where
some connected closed nilpotent subgroup of $G$ is transitive on $G/K$.

Example: $G/K$ where $G = H_m\rtimes K_m$ and $K = K_m$\,,
where $K_m$ occurs in the table above.

Fact: Let $G/K$ be commutative.  If a conn closed nilpotent 
subgroup $N$ of $G$ is transitive then $N$ is the nilradical
of $G$, $N$ is abelian or $2$--step nilpotent, and $G = N\rtimes K$.

In particular: Commutative nilmanifolds have form $G/K$ where
$G/K = (N\rtimes K)/K$,
$N$ is not so different from the Heisenberg group and $K \subset \Aut(N)$.

More examples: a commutative nilmanifold $(G = N\rtimes K)/K$ is
{\it irreducible} if $K$ acts irreducibly on $\gn/[\gn,\gn]$,
{\it maximal} if it is not of the form 
$(\widetilde{G}/\widetilde{Z},\widetilde{K}/\widetilde{Z})$ with
$\{1\} \ne \widetilde{Z} \subset \widetilde{K}$ central in $\widetilde{G}$.
They have been classified by Vinberg.  Let $\gn = \gz + \gw$ where $\gz$
is the center and $\Ad(K)\gw = \gw$.  Then $\gz$ is the center of $\gn$,
$\gw \cong \gn/[\gn,\gn]$ as $K$--module, and the classification is
\bigskip

{\footnotesize
\centerline{
\begin{tabular}{|r|c|c|c|}
\hline \hline
 & Group $K$ & $\gw$ & $\gz$   \\ \hline
1 & $SO(n)$ & $\R^n$ & $\Skew\R^{n\times n} = \gs\go(n)$   \\ \hline
2 & $Spin(7)$ & $\R^8 = \O$ & $\R^7 = \Im\O$    \\ \hline
3 & $G_2$ & $\R^7 = \Im\O$ & $\R^7 = \Im\O$   \\ \hline
4 & $U(1)\cdot SO(n), n \ne 4$ & $\C^n$ & $\Im\C$   \\ \hline
5 & $(U(1)\cdot) SU(n)$ & $\C^n$ & $\Lambda^2\C^n \oplus\Im\C$ \\ \hline
6 & $SU(n), n$ odd & $\C^n$ & $\Lambda^2\C^n$  \\ \hline
7 & $SU(n), n$ odd & $\C^n$ & $\Im\C$  \\ \hline
8 & $U(n)$ & $\C^n$ & $\Im \C^{n\times n} = \gu(n)$ \\ \hline
9 & $(U(1)\cdot) Sp(n)$ & $\H^n$ & $\Re \H^{n \times n}_0 \oplus \Im\H$   
        \\ \hline
10 & $U(n)$ & $S^2\C^n$ & $\R$  \\ \hline
11 & $(U(1)\cdot) SU(n), n \geqq 3$ & ${\Lambda}^2\C^n$ & $\R$ 
        \\ \hline
12 & $U(1)\cdot Spin(7)$ & $\C^8$ & $\R^7 \oplus \R$  \\ \hline
% \end{tabular}
% \begin{tabular}{|r|c|c|c|}
% \hline \hline
%  & Group $K$ & $\gw$ & $\gz$   \\ \hline
13 & $U(1)\cdot Spin(9)$ & $\C^{16}$ & $\R$  \\ \hline
14 & $(U(1)\cdot) Spin(10)$ & $\C^{16}$ & $\R$  \\ \hline
15 & $U(1)\cdot G_2$ & $\C^7$ & $\R$  \\ \hline
16 & $U(1)\cdot E_6$ & $\C^{27}$ & $\R$  \\ \hline
17 & $Sp(1)\times Sp(n), n \geqq 2$ & $\H^n$ & $\Im \H = \gs\gp(1)$  
        \\ \hline
18 & $Sp(2)\times Sp(n)$ & $\H^{2\times n}$ &
        $\Im \H^{2\times 2} = \gs\gp(2)$  \\ \hline
19 & $(U(1)\cdot) SU(m) \times SU(n)$ &  &   \\
   & $m,n \geqq 3$ & $\C^m\otimes \C^n$ & $\R$  \\ \hline
20 & $(U(1)\cdot) SU(2) \times SU(n)$ & $\C^2 \otimes \C^n$ &
        $\Im \C^{2\times 2} = \gu(2)$   \\ \hline
21 & $(U(1)\cdot) Sp(2) \times SU(n), n \geqq 3$ & $\H^2\otimes \C^n$ & $\R$ 
        \\ \hline
22 & $U(2)\times Sp(n)$ & $\C^2 \otimes \H^n$ & $\Im \C^{2\times 2} = \gu(2)$ 
        \\ \hline
23 & $U(3)\times Sp(n), n \geqq 2$ & $\C^3 \otimes \H^n$ & $\R$   
        \\ \hline
\end{tabular}
}}

\bigskip

\noindent
where the optional $U(1)$ is required 
in (5) when $n$ is odd, 
in (11) when $n$ is even,
in (19) when $m = n$, 
in (20) when $n = 2$, and
in (21) when $n \leqq 4$.
Here (9) was the first known case where $G/K$ is not weakly symmetric
(Lauret).

To make this explicit one needs to know the positive definite spherical
functions and the Plancherel measure.

In the connected
Lie group cases, the $K$--spherical functions on $G$ are just the 
joint eigenfunctions for $\cD(G,K)$, and in many cases this is how one finds
them.  We look at a few of those cases.

Case $(G,K) = (\R^n\rtimes K, K)$ where $K$ is transitive on the unit
sphere in $\R^n$.  Then the invariant differential
operators are the polynomials in the Laplacian 
$\Delta = -\sum \partial^2/\partial x_i^2$, and the 
$K$--spherical functions are the radial eigenfunctions of $\Delta$ for
real non-negative eigenvalue.  They are the
$$
\varphi_\xi(x) = (\Vert\xi\Vert r)^{-(n-2)/2}J_{(n-2)/2}(\Vert\xi\Vert r)
$$
where $r = \Vert x \Vert$ and $J_\nu$ is the Bessel function of first
kind and order $\nu$.

Case $(G,K) = (H_n\rtimes U(n),U(n))$.  In the coordinate $(z,w) \in \Im\C 
+ \C^n = H_n$ the invariant differential operators are the polynomials in
$\partial/\partial z$, $\square = -\sum \partial^2/\partial w_i^2$ and
$\overline{\square}$.  The positive definite spherical functions 
corresponding to $1$--dimensional representations are the
$$\varphi_\xi(z,w,k) = \tfrac{2^{n-1}(n-1)!}{(||\xi||\, ||w||)^{n-1}}
        J_{n-1}(||\xi||\, ||w||) \text{ for } 0 \ne \xi \in \gw^*$$
and those for infinite dimensional representations are the
$$\varphi_{\zeta,m}: (z,w,k) \mapsto
\begin{cases} e^{i\zeta(z)} L_m^{(n-1)}(\zeta(z)||w||^2)
                      e^{-\zeta(z)||v||^2/4} 
              \text{ if } \zeta(i) > 0, \\
              \phantom{x} \overline{\varphi_{-\zeta,m}(z,w,k)}
                      \phantom{XXXXXXXXi}\text{ if } \zeta(i) < 0,
\end{cases}$$
where $\zeta\in (\Im\C)^*$ and $L_m^{(n-1)}$ is the generalized Laguerre 
polynomial of order $n-1$ normalized to $L_m^{(n-1)}(0)=1$.

\bigskip
\bigskip

\begin{tabular}{l}
Department of Mathematics \\
University of California \\
Berkeley, California 94720--3840, U.S.A.\\
\\
{\tt jawolf@math.berkeley.edu}
\end{tabular}
\end{document}